\documentclass[11pt]{amsart}

\usepackage{amssymb, graphicx}
\usepackage{graphicx}

\usepackage[all]{xy}

\usepackage{amsfonts}
\usepackage{amsmath}
\usepackage{latexsym}
\usepackage{amscd}
\usepackage{verbatim}

\allowdisplaybreaks[1]

\newtheorem{theorem}{Theorem}[section]

\theoremstyle{remark}

\begin{document}
\title[What happened to PI theory?]{What happened to PI theory?}

\author{Claudio Procesi}

\address{C. Procesi, Dipartimento di Matematica, Sapienza
 Universit\` a di Roma,  Italy}  \email{procesi@mat.uniroma1.it}

\begin{abstract} This paper is a request for help to clarify a rather confusing state of the Theory.\medskip

\hskip4cm{\em dir\`o sia il peccato che il peccatore}\end{abstract}

\thanks{2010 {\em Math. Subj. Class.} Primary 16R60. Secondary 16R10, 16R30. }
\maketitle
\section{Introduction}
The notion of algebraic identity is quite general in Algebra, here we restrict to Polynomial Identities of associative algebras.   

Formally, given an algebra $R$ over a commutative ring $F$ a polynomial identity of $R$ is a polynomial $f(x_1,\ldots,x_m)$ in non commutative variables with coefficients in $F$ which vanishes identically when evaluated in $R$. 

 I will assume for simplicity that $F$ is an infinite field although it is quite interesting to consider other cases as for instance $F=\mathbb Z$ the integers.
\medskip

 I started my career in this field  under the guidance of Amitsur and Herstein and later left it for other areas. Recently I decided that I wanted to understand what had happened to PI theory since I left it several years ago  and found a somewhat sore state of the Theory,  many very sketchy papers,  sometimes truly obscure, some expositions in books containing gross mistakes and a general impression that some of the fundamental theorems presented have not been checked   by  anybody seriously.  
 
 Of course the problems of PI theory have not the same impact in Mathematics as the Riemann hypothesis or the Poincar\'e conjecture, so there is not a large community of strongly motivated scientists who aim to have a firm grasp on the theory, nevertheless PI theory, even if not so central, is a beautiful part of mathematics and deserves a better treatment that the one I have observed. Here I want to point out  where, in my opinion, the Theory needs a deep clarification in order to understand if results which have been celebrated and widely used are truly proved or even true. 
 
\section{A brief history}   Without pretending to give a historical account  which would take a much more serious effort and apologizing for all the omissions, one can say that the theory of Polynomial Identities of associative algebras which for short we call {\em PI theory}, has developed in the following steps.

 \subsubsection{ Kuros problem} The theory starts at the end of the 40's of last century with work of Jacobson,  \cite{jacobson1}, Kaplansky \cite{kap} and Levitzki \cite{Lev} achieving the goal of proving the bounded Kuros problem, that is 
 
 {\em A finitely generated  associative algebra $R$ over a field $F$ in which every element satisfies an algebraic equation of bounded degree is finite dimensional}.
 
 As the reader will realize this is an analogue of the bounded Burnside problem for groups which instead has a negative answer.\medskip
 
 The approach to the Kuros problem, which we may call the {\em American Israeli school}, is based on structure theory of non commutative rings, an important step in research at that time, summarized by the books of Albert and Jacobson.  One notices that an  algebra $R$ over a field $F$ in which every element satisfies an equation of bounded degree satisfies a polynomial identity and then a crucial Theorem of Kaplansky shows that a primitive PI algebra is a finite dimensional central simple algebra, this with work of Levitzki on the nil radical of a PI algebra gives enough information to conclude.\medskip
 
\subsubsection{The embedding problem}  A second problem, again due to a Russian algebraist Malcev  connected with the previous one was 
 
 {\em characterize associative rings which can be embedded in a finite matrix ring $M_n(A)$ over some commutative ring $A$.}
 
 Here  stands the  fundamental result of Amitsur--Levitzki  that the ring $M_n(A)$  satisfies the standard  polynomial $S_{2n}$.   For any integer $h$   the standard polynomial $S_{h}$ is the element of the free algebra  in the variables $x_1,\ldots,x_h$ given by the formula
$$S_h(x_1,\ldots,x_h):=\sum_{\sigma\in \mathfrak S_h}\epsilon_\sigma x_{\sigma(1)}\ldots x_{\sigma(h)} $$ where $\mathfrak S_h$ is the symmetric group of permutations on $\{1,\ldots,h\}$ and $\epsilon_\sigma$  denotes the sign of the permutation $\sigma$.\smallskip

In general thus a necessary condition  for an $F$ algebra to be embedded in a finite matrix ring $M_n(A)$ over some commutative $F$ algebra $A$   is that it should satisfy all polynomial identities of matrices   $M_n(F)$. The nature of the space of all polynomial identities is quite mysterious and not fully described except for $n=1,2$.

It was soon realized by Amitsur, Small and others (cf. \cite{irsmall}) that this condition is not sufficient, since there are other ring theoretical conditions (chain conditions) which must be satisfied; but in fact, as far as I am aware,  no complete necessary and sufficient condition is known.

The situation changes drastically if one takes the point of view of universal algebra, consider not only the structure of associative algebras but adds the axioms of a trace and works with {\em algebras with trace}. Then one has, in characteristic zero, that a necessary and sufficient condition for an algebra with trace $R$ to be embedded in a finite matrix ring $M_n(A)$  (compatible with the trace) is that $R$ should satisfy the $n^{th}$ Cayley--Hamilton identity (one can formulate possible analogues in all characteristics but at the moment no precise results are known). 

A development of this point of view leads to the {\em Theory of Cayley--Hamilton algebras}, cf.  \cite{P5}.
\medskip

\subsubsection{Geometry of representations} In fact the Theorem of Kaplansky  points out to the fact that PI theory is strictly connected to the study of finite dimensional representation of algebras, a theory which has a strong geometric flavor as shown by M. Artin and myself in the 60's and which can be considered as a second development of PI theory. 

Here one can mention the Nullstellensatz, that is the statement that, for a finitely generated PI algebra $R$, the ideal of elements vanishing in every irreducible finite dimensional representation coincides with the nil radical and the powerful theorem of M. Artin, \cite{Artin} and \cite{Sche},  which states that an algebra which  satisfies all polynomial identities of matrices   $M_n(F)$ and has no irreducible representation of dimension strictly smaller than $n$ is an Azumaya algebra.  Azumaya algebras are essentially {\em non split forms} of matrix algebras, as for instance finite dimensional simple algebras. In particular, such an algebra behaves in a very similar way to a matrix algebra  and can  be {\em split}; that is, it can be naturally embedded in a  matrix ring $M_n(A)$ over some commutative $F$ algebra $A$ (same $n$), the embedding should be faithfully flat as to preserve properties of the algebra.

Aside from the special case of Azumaya algebras which are very similar to commutative algebras  many difficult questions arise in general. With respect to the  Nullstellensatz,  many questions remained open at that time, in particular, on the more combinatorial side,  to understand if the nil radical of  a finitely generated PI algebra $R$ is nilpotent.
\medskip

I first mentioned the American Israeli school, to which one may oppose to some extent the Russian school, just one name here as PI theory is concerned, that of Shirshov, whose work was not known in the west due to  the cold war and major difficulties in communication between these two worlds at that time.  This school, with a strong combinatorial approach, was later able to unravel the mysteries of the Burnside problem, but this is another story and I leave it to more competent persons.

 A major breakthrough of combinatorial nature, but with a profound impact on structure theory, was the discovery, independently by Formanek \cite{formanek3} and Razmyslov \cite{Raz}, of central polynomials for matrices, that is non commutative polynomials which do not vanish identically on   $M_n(F)$  but always take scalar values.  It is quite remarkable that here we do not even know in general which are the central polynomials of minimum degree, something which should parallel the Amitsur--Levitzki theorem on the standard identity.
 
  \subsubsection{$T$--ideals}  When we consider the symbolic calculus associated to some algebra $R$, we are faced to study the free algebra $F\langle X\rangle$ in some  variables $X$  modulo the ideal $Id(R)$ of all polynomial identities of $R$.
  
  Ideals of type $Id(R)$ in  the free algebra $F\langle X\rangle$ are characterized by the property of being closed under the operation of substituting the variables $x_i$ with polynomials $g_i\in  F\langle X\rangle$. In  more formal language the semigroup of endomorphisms of  the free algebra $F\langle X\rangle$, is given by substituting the variables $x_i$ with polynomials $g_i\in  F\langle X\rangle$. We call  such semigroup $T_X$ or just $T$  and then call a $T$ ideal an ideal of $F\langle X\rangle$  stable under the action of $T$.  The quotient algebra $F\langle X\rangle/I$ modulo a $T$ ideal is then also a free algebra, but not for the category of all algebras but for that of algebras satisfying the identities of $I$. This is just like what happens for commutative algebras where the free algebras are the algebras of commutative polynomials.\smallskip
  
  Usually in order to avoid confusion  the quotient algebra $F\langle X\rangle/I$ modulo a $T$ ideal is called a {\em relatively   free algebra}.  Here then several major difficulties of PI  theory appear, from a purely algorithmic and combinatorial point of view, while the free algebra on one hand, and the polynomial algebra on the other have a clear  description in term of a basis of {\em monomials}, this is not the case for relatively free algebras.  Only for very special examples such a description is possible and even the dimension of the single graded components of a relatively free algebra is usually not explicitly computable.
  \smallskip
  
  The second major difficulty  is in that the analogue of {\em Hilbert basis theorem} does not hold, in other words  relatively   free algebras are neither right nor left Noetherian. Therefore, of special interest, is a conjecture of Specht that 
  
  {\em A $T$ ideal is finitely generated when considered not only as a bimodule over $F\langle X\rangle/I$ but also under the action of the semigroup $T$}.
  
  A solution to this problem has been presented by Kemer \cite{kemer} as a consequence of a deeper analysis of the structure of $T$ ideals, this is a crucial part of the Theory whose presentation I find unsatisfactory, to say the least, and on this I will concentrate my discussion in the next section.
  \smallskip
  
  In the theory of $T$ ideals a special role plays the $T$ ideal of  polynomial identities of matrices   $M_n(F)$. I already mentioned that even this ideal is not fully described except for $n=1,2$, nevertheless many beautiful properties of the corresponding relatively free algebra are known. The Theory was initiated by Amitsur who showed that this algebra is an integral domain with a quotient division algebra of rank $n^2$ over its center. In my Ph. D. thesis I started to investigate this object as {\em  the algebra of generic matrices} and from the point of view of invariant Theory and started to look at the larger algebra in which coefficients of the characteristic polynomial are added,  this later led to the Theory of Cayley--Hamilton algebras. Let us for future reference denote by $R_{n,m}$ the algebra of $m$ generic $n\times n$ matrices, it is embedded naturally in the matrix algebra $M_n(A)$ where $A$ is the polynomial ring of functions on $m$--tupels of $n\times n$ matrices.
  
 Denote by $T_{n,m}$ the commutative algebra generated by all the coefficients of characteristic polynomials of elements of $R_{m,n}$. One can prove that  $T_{n,m}$ is a finitely generated commutative algebra and $R_{m,n}T_{n,m}$ is a finitely generated module over the commutative algebra $T_{n,m}$. In particular both  $T_{n,m}$ and $R_{m,n}T_{n,m}$ are Noetherian (contrary to $R_{m,n} $), cf. \cite{procesi2}.

   \subsubsection{Growth}  I already mentioned that the dimension of the single graded components of a relatively free algebra is usually not explicitly computable. Nevertheless some strong and useful results can be obtained if one limits the analysis to  asymptotic statements, that is estimates on the growth.   Such results have been obtained by several authors  and depend on the representation Theory of the symmetric group. An important role plays the theorem of Regev \cite{regev1}, that the tensor product of two PI algebras satisfies a PI, and the study of growth gives rise to useful invariants of PI algebras as shown by Giambruno--Zaicev. These investigations have been quite active in the last years and I will not discuss them since I have no concern on their  presentation, in fact one has a good  reference in the book of  Giambruno--Zaicev, cf. \cite{giambruno}.
  
  \subsubsection{Trace identities} I already mentioned that the introduction of traces  can change drastically the theory, of course in general a ring has no natural trace and, although one can add formally a trace, the resulting object has no natural identity involving this formal trace. 
  
  Introducing {\em algebras with trace}  one has a natural notion of free algebras and trace identities and then the first major result, discovered independently by Razmyslov and myself \cite{procesi}, \cite{Raz}, which can be thought of as the {\em second fundamental theorem for matrix invariants}   is that all trace identities of $n\times n$ matrices can be deduced, in characteristic 0, from the Cayley--Hamilton identity.
  
  The algebra $R_{m,n}T_{n,m}$ is the relatively free algebra, in the category of algebras with trace satisfying the identities of $n\times n$ matrices, in $m$ variables.\smallskip
  
  I started from this result in order to build the Theory of Cayley--Hamilton algebras while Razmyslov based on this his solution of the problem on the nil radical, that is proving, in characteristic 0 that the radical of a finitely generated algebra is nilpotent.  This I discuss in the next paragraph.
  
  \subsubsection{The Razmyslov--Braun--Kemer Theorem}
 The question whether the radical of a finitely generated algebra is nilpotent arises naturally from the Nullstellensatz and from the Levitzki--Amitsur  theory on nil ideals in PI rings. They proved that, if $R$ satisfies a proper multilinear identity of degree $d$, then if $T\subset R$ is a nil subring we have $T^{[d/2]}\subset N(R)$ where $N(R)$ is the sum of all nilpotent ideals.
  
  In a fundamental paper Razmyslov explains how to use trace identities to overcome the problem that finitely generated PI algebras are not noetherian.  For this he needs to have some control on the identities and so introduces what he calls the {\em Capelli identity}
  \begin{equation}
\label{Cap}C_m:=\sum_{\sigma\in S_m}y_0x_{\sigma(1)}y_1x_{\sigma(2)}y_2\ldots y_{m_1}x_{\sigma(m)}y_m
\end{equation}
I find this name somewhat misleading since in the classical literature the Capelli identity is a certain identity between differential operators. Indeed the classical Capelli identity and the identity introduced by Razmyslov play a similar role as they allow to reduce the structure of a $T$ ideal modulo such an identity to the structure of its intersection with the free algebra in $m-1$ variables, something that, in classical invariant Theory is called the {\em method of primary covariants} and, in modern representation Theory, the {\em Theory of highest weight vectors}.

A major remark of Razmyslov has been to show that if one considers the algebra of generic matrices $R_{m,n}\subset  R_{m,n}T_{n,m}$ this inclusion has a {\em non zero conductor} that is a non zero ideal $I$ of both  algebras.

This implies that $R_{m,n}$ contains an ideal $I$ which is a finitely generated module over a finitely generated commutative algebra  but in fact, if $\sqrt{I}$ denotes the radical of $I$,  one has that $R_{m,n}/\sqrt{I}=R_{m,n-1}.$

This important observation is the key to make inductions and allows Razmyslov to prove that 

{\em In characteristic 0   the radical of a finitely generated algebra satisfying some Capelli identity is nilpotent. }

 Unfortunately Razmyslov's paper is very compressed, and sketchy at points, so that it has remained for a long time an obstacle to  many people who took it by faith (I am preparing some extended presentation of this and other papers on the subject which are even more compressed).

The argument of Razmyslov needs the assumption that  the algebra   satisfies some Capelli identity and hints at a possible argument in positive characteristic.

The first problem was solved by Kemer who showed that,  any finitely generated algebra  satisfies some Capelli identity. This Theorem may be deduced, in characteristic 0,  from growth estimates.

The Theorem in a characteristic free way was later  proved by Braun.

The theory is complemented by the Theorems of Lewin, which gives a nice construction of an algebra  having as $T$ ideal of identities the product of two $T$--ideals.  This follows the steps of Fox calculus and the theory of Magnus for metaabelian groups.

Then the results of Giambruno--Zaicev, which tell us that  if $\Gamma_i$  is the $T$ ideal of identities of $i\times i$ matrices, then a product $\Gamma_{i_1}\Gamma_{i_2}\ldots \Gamma_{i_m}$  of these $T$ ideals is the $T$ ideal of idenitities of {\em block triangular}  matrices, with blocks of sizes $i_1,i_2,\ldots,i_m$  on the diagonal.
\subsection{ The Theory of Kemer}

In the 80s Kemer  introduced a  large set of new ideas in the Theory  and presented some major results.

For one thing,  by a systematic use of the infinite Grassmann algebra and the idea of Grassmann envelope of a superlagebra, he was able to treat on a  parallel footing, symmetry and antisymmetry  in the Theory of cocharacters, thus treating algebras with PI which do not necessarily satisfy a Capelli identity. On the other hand, he has attacked the Specht problem through the question of representability of relatively free algebras.

\subsubsection{Representable algebras}   I have already mentioned that one of the first inputs of PI theory is the question to characterize algebras which can be embedded in  some algebra $M_n(T)$  of matrices over a commutative ring $T$, one says if this is possible that $R$ is representable.  In particular this question has been asked for an algebra $R=A[a_1,\ldots,a_k]$ finitely generated over some commutative ring.  There are variations of this question, one can ask if  $R$   can be embedded in  some algebra $M_n(K)$  of matrices over a field $K$. It is then easy to see that if  $R=F[a_1,\ldots,a_k]$  and $F   $ is also a field  the two questions are equivalent, that is if  $R$ can be embedded in  some algebra $M_n(T)$  of matrices over a  commutative ring then it   can also be embedded in  some algebra $M_m(K)$  of matrices over a field $K$ (usually $m>n$).

We have already noticed that a   finitely generated PI algebra $R=A[a_1,\ldots,a_k]$ need not be representable, since there are some necessary chain conditions on annihilators to be satisfied that are not satisfied in general.

Thus, the following result established by Kemer is quite surprising.  Let us consider algebras over a field $F$.
\begin{theorem}\label{}
A finitely generated relatively free algebra is representable.
\end{theorem} This result can be stated also in a different but equivalent form, 
\begin{theorem}\label{}
Given any finitely generated PI algebra $R$ over an infinite field $F$ there is a finite dimensional algebra over $F$ which satisfies the same polynomial identities as $R$.\end{theorem}
When two algebras satisfy the same polynomial identities one says that they are {\em PI equivalent}, thus the statement is that a finitely generated PI algebra is PI equivalent to a finite dimensional algebra.\smallskip

In fact Kemer presents the result already in the generality of superalgebras so to deduce even a more general statement that any PI algebra $R$ (not necessarily finitely generated) is PI equivalent to the Grassmann envelope of a finite dimensional superalgebra.

These results are then instrumental for his solution of the Specht problem.\smallskip
  
  Unfortunately Kemer's monograph is written in such a sloppy form that in reading  it I found a lot of statements which I could not understand or which at first sight sound like nonsense, so I certainly am  unable  to say if his presentation is   correct or complete, but even to claim that he is wrong is quite hard  (as people say of some Theory, {\em not even wrong}).
  
  To be very precise, I stopped reading when I reached lemma 2.4, a technical statement for the crucial lemma 2.5.
  After many attempts to decode lemma 2.4 I gave up.\smallskip

Many people have written  papers based on this work or generalizing it, but my impression has been that in all cases the work of Kemer was given for granted and the papers were built on it.

Asking around no one I know seems to have read the proof  that a finitely generated PI algebra is PI equivalent to a finite dimensional algebra. 

To make things worse there is the book of Belov--Rowen \cite{belov} in which they present a proof of this result, unfortunately at some point they use the wrong formula on vector spaces distributing sums and intersection, but even if this can be corrected (and I think it can) the rest of the proof behaves like some of these  rivers of southwestern Asia and northern Africa, the {\em wadi}  which at some point just disappear in the sand.

One thing which makes me particularly worry is that I have the feeling that, the place in which Belov--Rowen get lost,  looks to me very suspiciously similar to the place in which I got lost trying to understand Kemer.

To be precise  the strategy seems to be that, given a  $T$ ideal  $\Gamma$ one tries to approximate $\Gamma$ from {\em below} as much as possible, thus finding an ideal $Id(A)\subset \Gamma$  of identities of a finite dimensional algebra, which shares with $\Gamma$ sufficiently many identities and non identities. This is measured by certain invariants which in   Belov--Rowen are called Kemer index, the proof is thus by induction on this index. The first step is thus to find  this $A$ so that $Id(A)\subset \Gamma$ have the same Kemer index, I understand  (in my own way)  this part.  What is quite unclear is how  you conclude.  Here the step seems to be a combination of a simple induction by lowering the Kemer index of  $\Gamma$ with a technique of {\em eliminating some identities}, this last part I just do not understand.\smallskip
 
  I am totally aware that, not understanding a paper, is just {\bf a fact} and it came to my mind a supposed episode, maybe just a metropolitan legend, at the Institute of Advanced Study in Princeton.
  
  The way in which the story has been told to me is the following:
  
  After a very difficult lecture of Dirac,  he asked as usual "Are there any questions?"
  The lecture was attended by the cream of $20^{th}$ century scientists, Herman Weyl, Albert Einstein, Kurt G\"odel. The story goes, Einstein said "I do not understand ... so and so."  Dirac replied  "This is not a question, it is a fact".
  
  \smallskip
  
  But, as far as Kemer's results are concerned, the question is not if Claudio Procesi understands this theorem or not;  but, the fact is that the PI community is very small and, when asking around half of the people, they confess to have not read any proof and the other half refers to a  book with a major mistake, the situation starts to be worrisome. 
  
  The situation worsens when trying to read some developments of the Theory.
  
  Undoubtably a major goal of the Theory is the solution of the Specht problem \cite{specht}, and in this respect  some generalizations have been presented by Grishin and Shchigolev (cf. \cite{gris}, \cite{Shig}) who claim  that not only every $T$ ideal of a finitely generated free algebra is finitely generated (as $T$ ideal) but in fact any $T$ stable submodule is also finitely generated. Recall that $T$ is the semigroup of endomorphisms  of the free algebra. So, the results look interesting and next I tried to read these papers.

 Now, in one of Grishin's papers, used by Shchigolev,  he makes the following definitions:

$k$ a field $k_r$ the ring of $r\times r$ matrices:

  he calls an algebra $A$ with 1 a " $k_r$ algebra" if  $A$ contains  $k_r$  as algebra with 1.

Somehow he does not seem to realize that this implies that $A=M_r(B)$  is also  $r\times r$ matrices over some ring $B$.

Next he defines a variety of algebras "$r$--closed,"    if it is generated by a $k_r$ algebra. 

This certainly is an interesting notion but then he makes the following statement which I find absurd,
I quote

"{\em It is easy to show that any $r$ closed variety is generated by its finite dimensional algebras with unity which are of the form $k_r+N$, where $N$ is the radical of  the algebra $k_r+N$}"  end quote.

For me an example of $k_r$ algebra is clearly $k_{rs}=k_r\otimes k_s$ for any $s$  so in whatever way you try to read this statement it looks to me as nonsense. Or am I missing something?

After two pages Grishin discusses exactly this case so I am in total confusion,  I have  tried to continue reading these papers but found them unreadable.

One hypothesis is that, since most, if not all, of these papers were originally written in Russian, it is possible that something is {\em lost in translation}, unfortunately my knowledge of Russian is probably reduced to at most 5 words and even if I am quite interested I am afraid that at my age I cannot engage in learning Russian just for the purpose of reading these papers.

On the other hand, I  am willing to bet that not much is lost in translation. A few years ago I wanted to read some results in algebraic topology  which were available only in Russian, google translator did not exist at the time, so I took a dictionary and, with some pain, a word here and there, a few formulas, I was able to understand the paper. Certainly my ability of translating from Russian is lower than that of even the worst  translator.\smallskip

I am quite afraid that, since I find these papers so irritating, I lower my  professional defenses and do not understand what is truly meant by these authors, but maybe they are just wrong.

But the worst possible situation is if they are {\bf not}  wrong, but nobody understands their work.  This is a sure way in which some mathematics may be in the future lost forever.

\subsubsection{A request}
Many of the difficulties of PI theory stem from the fact that, in particular in the more combinatorial aspects, there are few general results that one can readily use, there are many involved special computations and identities in which one can get easily lost.

I think a first effort should be to lay down the foundations of the combinatorial part in a more systematic way.

One should as much as possible try to clarify which identities can be deduced from some basic ones.
 The problem seems to me to be that $T$ ideals are a much more complicated object than representations,  in the tensor algebra one has a nice decomposition into isotypic components of the linear group but to understand the contributions to a $T$ ideal is quite difficult.
 On the other hand it cannot be that reading a paper on PI is like trying to decipher hieroglyphics, to say the least because the interest for humanity is not at the same level.
 
 A suggestion that I received from Lance Small and which seems quite interesting is to open a blog dedicated to a discussion on PI Theory, I hope  many people may be interested in this discussion.

\printindex

\end{document}